\documentclass{amsart}

\usepackage{amsfonts}

\makeatletter
\renewcommand{\@evenhead}{\hfil \thepage}
\renewcommand{\@oddhead}{\hfil \thepage}
\renewcommand{\@evenfoot}{}
\renewcommand{\@oddfoot}{}
\makeatother
\theoremstyle{plain}

\begin{document}
\title[Controlling alliances through executing pressure]{Controlling
alliances through executing pressure}
\thanks{$^{\ast }$Corresponding author}
\author{L. A. Khodarinova$^{1}$}
\thanks{$^{1}$SPMMRC, School of Physics and Astronomy, University of
Nottingham, Nottingham NG7 2RD UK}
\author{J. M. Binner$^{2\ast }$}
\thanks{$^{2}$Dept. of Strategic Management, Aston Business School, Aston
University, Birmingham B4 7ET UK}
\author{L. R. Fletcher$^{3}$}
\thanks{$^{3}$School of Computing and Mathematical Sciences, Liverpool John
Moores University, Liverpool L2 3AA UK}
\author{V. N. Kolokoltsov$^{4}$}
\thanks{$^{4}$School of Computing and Mathematics, Nottingham Trent
University, Nottingham NG1 4BU UK}
\author{P. Whysall$^{5}$}
\thanks{$^{5}$Nottingham Business School, Nottingham Trent University,
Nottingham NG1 4BU UK}

\begin{abstract}
In this paper the standard prisoners' dilemma is embedded in environmental
conditions in which the interaction takes place. This provides a theoretical
background to the analysis of the empirical studies which indicate that
including additional factors when considering an alliance is very important.
We show that such an approach, though theoretically simple, provides a
powerful tool for suggesting successful strategies for forming, maintaining
and withstanding a rival attack on an alliance, projecting an alliance
possible success and determining which measures must be taken to maintain an
alliance in a changing environment. It also provides insight into the
possibility of preventing alliances in rivals. The relative simplicity of
this framework means that the approach can be easily applied when the real
life decisions are taken.
\end{abstract}

\maketitle

\section{Introduction}

Recently there has been an increasing interest in understanding the
instruments responsible for success in business alliances. In a highly
competitive and globalized market alliances and joint ventures are becoming
increasingly popular by allowing companies to share risks, develop new
markets, acquire technological advances and so on. There has also been a
significant increase in alliances between competitors$^{1,2}$ as such
ventures not only improve their market performance but also reduce costs due
to eliminating rivalry hazards. Unfortunately, empirical research$^{2-8}$
shows that business alliances have high failure rates (at or higher than
50\%). The reason for such failure has a dual nature$^{2,9}.$ On one side
non opportunistic reasons such as managerial complexity and coordinating
costs play a significant role. The other important reason is interfirm
rivalry as each partner is acting to maximize its own benefits. This paper
is aimed at addressing the second reason.

In order to explain the cooperative behaviour when the risk of opportunistic
behavior is high game-theoretic descriptions of social interactions are
often\ used$^{10,11,12}.$ The single-interaction prisoners' dilemma has been
widely used$^{11,12}$ as a generic model of such\ interactions. Analysis of
the single interaction prisoners' dilemma leads to the conclusion that
cooperative behaviour cannot be rational; nevertheless it is observed by
many empirical studies. The main approach to resolving this paradox has been
based on considering the iterated prisoners' dilemma, i.e. considering the
game in the long term context$^{9,10,11}.$ This changes the structure of the
game and cooperation becomes rational$^{13,14}.$ The basic idea is that by
including future interactions \textquotedblleft punishing
strategies\textquotedblright\ may be used, i.e., \textquotedblleft if you do
not cooperate now, you will be punished in the future\textquotedblright .
Unfortunately, the empirical studies indicate$^{10}$ that companies using
such strategies are not very successful in maintaining alliances. Nor are
such strategies evolutionary stable and if a player operates in a
noncooperative society no alliance can evolve in the long run$^{15}.$

From a game-theoretic point of view there exist many different approaches to
overcoming these difficulties$^{16-24}.$ Such approaches, most of which are
based on computer simulations$^{17-23},$ in particular include considering
complex strategies in an evolutionary model$^{17,18}.$ For example players
may be allowed to learn from experience$^{19}.$ It is shown$^{19}$ that
mutual cooperation can be maintained when players have a primitive learning
ability. It was shown that under proposed learning evolution some
cooperative strategies can invade not only unconditional cooperation, Tit
for Tat and Pavlov strategies but also noncooperative strategies. Another
approach is to embed the prisoners' dilemma into a spatial context$^{20,21}.$
For example$^{20},$ the version of the iterated prisoners' dilemma with only
unconditionally cooperating and unconditionally defecting players
interacting with the immediate neighbours was considered. It was shown$^{20}$
that such a model can generate chaotically changing spatial patterns, in
which cooperators and defectors both persist indefinitely. The pattern
generated by groups of cooperators exhibits$^{21}$ the scale-invariance
which is typical of self-organized criticality$^{22}$. Multi players games
are also considered$^{23}$ with solutions found by computer simulations.
There is also a direction of research based on the idea of private
information available to the players. That is, each player may have some
information about themselves or state of the game which is not available to
other players. The recent developments$^{24}$ in this area have revealed the
possibility of cooperation under condition of private information. Providing
an important framework for understanding evolution of cooperation in Nature,
such approaches however are not easily applicable to the situations when the
real life business alliances must be managed. Let us notice nevertheless
that all of these studies indicate that cooperative behaviour becomes
understandable when the interaction is considered in a wider context (such
as length of interaction, for example).

Surprisingly, there are few theoretical studies that suggest modeling
cooperation by the games other than prisoners' dilemma. Empirical research
on the other side shows that including additional factors when considering
an alliance is very important and may suggest successful strategies for
forming and maintaining an alliance. For example$^{1},$ it is shown that
engaging in unilateral commitments significantly increases an alliance
success rate. It is argued$^{1}$ that unilateral commitments change the
payoff structure of the game in such a way that mutual defection is no
longer the Nash Equilibrium in the single interaction game. Unfortunately
the lack of general game theoretic description of such and similar
situations leaves it unexplained why some unilateral commitments lead to an
alliance in one situation but not in others.

In this paper we try to overcome this deficiency by embedding the standard
Prisoners Dilemma in environmental conditions in which the interaction takes
place. We show that such an approach, thought theoretically simple, provides
a powerful tool for suggesting successful strategies for forming,
maintaining and withstanding a rival attack on an alliance, projecting an
alliance's possible success and determining which measures must be taken to
maintain an alliance in a changing environment. It also provides insights
into the possibility of preventing alliances in rivals. The relative
simplicity of this framework means that the approach can be easily applied
when the real life decisions are taken.

The main idea of our approach is to evaluate the external conditions under
which the alliance is forming. We suggest that the majority of such
conditions can be evaluated as either unchangeable or changeable pressure
(benefits) to the players. As far as the unchangeable factors are concerned
they can be considered as the base payoffs of the single interaction game
(for example prisoners' dilemma in our case). We model changeable conditions
by introducing a \textquotedblleft third player\textquotedblright . The
\textquotedblleft third player\textquotedblright\ can be irrational$^{25}$
and, for example, can represent \textquotedblleft Nature\textquotedblright ,
that is environmental or other conditions that have an impact on the game
but do not depend on the behaviour of the first two players (for example, a
highly competitive market or a tax regime imposed by a government). This
could also include pressure that a company decides to impose on itself by
engaging in unilateral commitments. The \textquotedblleft third
player\textquotedblright\ can also be rational. In this context there might
be a real third player who interacts with the first and the second players
(engaged in a prisoners' dilemma type game for example), possibly to a
degree he can choose.

This model correlates well with the empirical studies showing that by
considering the two players not in isolation but in the context of their
interaction with the \textquotedblleft third player\textquotedblright\ a
cooperative Nash Equilibrium between the first two players may emerge. We
investigate the effect which the existence of the \textquotedblleft third
player\textquotedblright\ has on the possibility of cooperation. This
provides mechanisms for controlling and managing the alliances via changing
the external conditions. Considering the \textquotedblleft third
player\textquotedblright\ as another company attacking the alliance we show
how to find an optimal strategy for the attack and also how such an attack
could be resisted and the alliance could be saved by good management.

\section{Modeling interaction under external pressure}

In this section we give an example of modelling an interaction in the
context of external conditions. As our example we consider two rival
companies engaged in the prisoners' dilemma game who, at the same time, are
withstanding an attack from a third player. In the model below we choose
certain unchangeable factors such as the value of the market or corporate
benefits. Clearly, it is done solely for the purpose of illustrating the
approach and the parameters appropriate to the situation under consideration
must be taken by the practitioners using this framework. For simplicity we
also consider a symmetric base interaction which also can be modified
without much difficulty. It is however important to note that the changeable
pressure is modelled continuously. Although in the current situation the
specific values of such pressure might apply, it is necessary to estimate
the impact of the pressure over the whole range of possible outcomes.

We begin by modelling an interaction with the prisoners' dilemma game in
which the players simultaneously choose an action from a set of two actions:
\textquotedblleft Ally\textquotedblright\ and \textquotedblleft
Fight\textquotedblright . Both players choosing \textquotedblleft
Ally\textquotedblright\ results in forming an alliance in which the players
share a market (resource) equally and obtain payoff of $l.$ Here $2l$
represents the value of the market (resource) share available. One player
choosing \textquotedblleft Fight\textquotedblright\ and the other
\textquotedblleft Ally\textquotedblright\ will be taken as the former player
is attacking the latter who does not fight back. Here $v$ represents the
value of the market (resource) share transferred as a result of the attack
from the \textquotedblleft allying\textquotedblright\ player to the
\textquotedblleft fighting\textquotedblright\ player. The parameter $x$
represents extra costs of running the business outside the alliance compared
with the cost of running the business in alliance (for example, corporate
benefits) If both players chose \textquotedblleft Fight\textquotedblright\
then each player is supposed to have an equal probability of obtaining the
whole market (resource). Here $f$ represents the cost of fighting. In order
for this game to be identified as a prisoners' dilemma we must have $v>x$
and $v>f$ and we also assume that $x\geq 0$ and $f\geq 0.$ The payoffs for
this game are summarized in the table below.%
\begin{equation}
\begin{tabular}{ccc}
$^{P_{1}\,}\backslash _{\,P_{2}}$ & Ally & Fight \\ 
Ally & $^{l\,}\backslash _{\,l}$ & $^{l-v-x\,}\backslash _{\,l+v-x}$ \\ 
Fight & $^{l+v-x\,}\backslash _{\,l-v-x}$ & $^{l-x-f\,}\backslash _{\,l-x-f}$%
\end{tabular}
\label{3Players_PD}
\end{equation}%
Building on the model let us now assume that there is a third player who
attacks the $i^{th}$ player with level $q_{i},$ $i=1,2.$ For the convenience
of the analysis we normalize the levels of the attack to be between zero and
one: $0\leq q_{i}\leq 1.$ If the third player attacks with level $q_{i}$ a
market share proportional to the appropriate level of attack is lost by the $%
i^{th}$ player. We assume that there may be a cost of being attacked which
is proportional to the level of the attack and will be expressed as $q_{i}c,$
where $c\geq 0.$ We also assume that withstanding an attack in an alliance
may be easier for a company than on its own, so that there maybe an extra
cost $q_{i}y,$ $y\geq 0,$\ of an attack if companies are not in an alliance.
Finally, we suppose that if players are in an alliance they share the cost
of the attack equally. These assumptions result in payoffs summarized in the
bi-matrix below.%
\begin{equation}
\begin{tabular}{ccc}
$^{P_{1}\,}\backslash _{\,P_{2}}$ & Ally & Fight \\ 
Ally & $^{l-lq_{1}-\frac{c}{2}\left( q_{1}+q_{2}\right) \,}\backslash
_{\,l-lq_{2}-\frac{c}{2}\left( q_{1}+q_{2}\right) }$ & $^{l+v-x-q_{1}\left(
l-v+c+y\right) \,}\backslash _{\,l-v-x-q_{2}\left( l+v+c+y\right) }$ \\ 
Fight & $^{l+v-x-q_{1}\left( l+v+c+y\right) \,}\backslash
_{\,l-v-x-q_{2}\left( l-v+c+y\right) }$ & $^{l-x-f-q_{1}\left( l+c+y\right)
\,}\backslash _{\,l-x-f-q_{2}\left( l+c+y\right) }$%
\end{tabular}
\label{3Players_new_payoffs}
\end{equation}%
It is fairly simple to obtain the pure strategy Nash Equilibrium conditions,
which are summarized in the Table I below.%
\begin{equation}
\begin{array}{c}
\text{Table I : Nash Equilibrium conditions for two players.} \\ 
\begin{tabular}{ccccc}
\hline
NE $^{P_{1}}\backslash _{P_{2}}$ & $^{Ally}\backslash _{Ally}$ & $%
^{Fight}\backslash _{Fight}$ & $^{Ally}\backslash _{Fight}$ & $%
^{Fight}\backslash _{Ally}$ \\ \hline
$_{Condition}$ & $%
\begin{array}{c}
q_{1}\geq \frac{2v-2x+cq_{2}}{2v+2y+c} \\ 
q_{2}\geq \frac{2v-2x+q_{1}c}{2v+2y+c}%
\end{array}%
$ & $%
\begin{array}{c}
q_{1}\leq 1-\frac{f}{v} \\ 
q_{2}\leq 1-\frac{f}{v}%
\end{array}%
$ & $%
\begin{array}{c}
q_{1}\geq 1-\frac{f}{v} \\ 
q_{2}\leq \frac{2v-2x+q_{1}c}{2v+2y+c}%
\end{array}%
$ & $%
\begin{array}{c}
q_{1}\leq \frac{2v-2x+cq_{2}}{2v+2y+c} \\ 
q_{2}\geq 1-\frac{f}{v}%
\end{array}%
$ \\ \hline
\end{tabular}%
\end{array}
\label{NE}
\end{equation}%
This indicates that the cooperative Nash Equilibrium can be obtained if the
levels of attack $q_{i}$ on both players is quite high. If $q_{i}$ are not
high enough the standard non cooperative Nash Equilibrium is observed. If $%
q_{i}$ vary significantly the weakly attacked player fights while the
strongly attacked player does not fight back. There could also\ be a mixed
strategy Nash Equilibrium as we will show in the example below.

As a next step we consider the rational third player and determine the
optimal levels of the attack $q_{i}.$ In this model we assume that before
the game~(\ref{3Players_new_payoffs}) is played the third player chooses $%
q_{i}$ or, equivalently, that the first two players are certain about the
strength with which they are attacked. This approach simplifies the analysis
of the game.\ The first and the second player make their choices at the
game~(\ref{3Players_new_payoffs}) which determine their own payoffs and also
the payoff obtained by the third player. The payoffs to the third player
will be given in Table II.%
\begin{equation}
\begin{array}{c}
\text{Table II: Payoffs to the third player.} \\ 
\begin{tabular}{ccccc}
\hline
NE $^{P_{1}}\backslash _{P_{2}}$ & $^{Ally}\backslash _{Ally}$ & $%
^{Fight}\backslash _{Fight}$ & $^{Ally}\backslash _{Fight}$ & $%
^{Fight}\backslash _{Ally}$ \\ \hline
$^{Payoff}$ & $^{\left( l-\sigma -\varepsilon \right) \left(
q_{1}+q_{2}\right) }$ & $^{\left( l-\sigma -\epsilon \right) \left(
q_{1}+q_{2}\right) }$ & $^{\left( l-\sigma \right) \left( q_{1}+q_{2}\right)
+v\left( q_{2}-q_{1}\right) }$ & $^{\left( l-\sigma \right) \left(
q_{1}+q_{2}\right) +v\left( q_{1}-q_{2}\right) }$ \\ \hline
\end{tabular}%
\end{array}
\label{3Players_payoffs_to_the_third}
\end{equation}%
Here $\sigma $ represents the cost of the attack, $\varepsilon $ represents
the additional cost due to attacking an alliance and $\epsilon $ represents
the decrement in the cost of the attack if it is made on the companies which
are fighting with each other.

We now consider examples showing how the optimal level of attack can be
chosen.

\section{Example}

Consider the game with the following parameters $v=5,$ $x=2,$ $y=2,$ $f=2$
and $c=3.$ Using formulae~(\ref{NE}) we obtain the following conditions for
the various Nash Equilibria.%
\begin{equation}
\begin{array}{c}
\text{Table III : Nash Equilibrium conditions for two players when} \\ 
v=5,x=2,y=2,f=2\text{ and }c=3. \\ 
\begin{tabular}{ccccc}
\hline
NE $^{P_{1}}\backslash _{P_{2}}$ & $^{Ally}\backslash _{Ally}$ & $%
^{Fight}\backslash _{Fight}$ & $^{Ally}\backslash _{Fight}$ & $%
^{Fight}\backslash _{Ally}$ \\ \hline
$_{Condition}$ & $%
\begin{array}{c}
q_{1}\geq \frac{3q_{2}+6}{17} \\ 
q_{2}\geq \frac{3q_{1}+6}{17}%
\end{array}%
$ & $%
\begin{array}{c}
q_{1}\leq \frac{3}{5} \\ 
q_{2}\leq \frac{3}{5}%
\end{array}%
$ & $%
\begin{array}{c}
q_{1}\geq \frac{3}{5} \\ 
q_{2}\leq \frac{3q_{1}+6}{17}%
\end{array}%
$ & $%
\begin{array}{c}
q_{1}\leq \frac{3q_{2}+6}{17} \\ 
q_{2}\geq \frac{3}{5}%
\end{array}%
$ \\ \hline
\end{tabular}%
\end{array}
\label{Ex_NE}
\end{equation}%
We can also find that if $\frac{3q+6}{17}\leq q_{1}\leq \frac{3}{5}$ and $%
\frac{3p+6}{17}\leq q_{2}\leq \frac{3}{5}$\ then there is a mixed strategy
Nash Equilibrium in which the first player and the second player choose
\textquotedblleft Ally\textquotedblright\ with probabilities $\frac{6-10q_{2}%
}{7q_{2}-3q_{1}}$ and $\frac{6-10q_{1}}{7q_{1}-3q_{2}},$ correspondingly.
For each combination of the parameters we have a unique Nash Equilibrium
solution except for the region where $\frac{3q+6}{17}\leq q_{1}\leq \frac{3}{%
5}$ and $\frac{3p+6}{17}\leq q_{2}\leq \frac{3}{5}$ for which there are
three Nash Equilibria: \textquotedblleft Ally\textquotedblright
-\textquotedblleft Ally\textquotedblright , \textquotedblleft
Fight\textquotedblright -\textquotedblleft Ally\textquotedblright\ and the
mixed Nash Equilibrium. It can be shown that in this region
\textquotedblleft Ally\textquotedblright -\textquotedblleft
Ally\textquotedblright\ is Pareto efficient meaning that both players obtain
their highest payoffs if the \textquotedblleft Ally\textquotedblright
-\textquotedblleft Ally\textquotedblright\ is played. The Nash Equilibrium
regions are shown in Figure 1.

Let us now fix parameters determining the payoff of the third player at $%
\varepsilon =1,\epsilon =1,\sigma =3$ and consider two scenarios: of high $%
(l=8)$ and low $(l=5)$ value of the market (resource). The payoff to the
third player depends on the strategies of the first two players~(\ref%
{3Players_payoffs_to_the_third}). If the Nash Equilibrium is not unique for
some range of the parameters $q_{i}$ we assume that the Pareto efficient
Nash Equilibrium is played. The payoffs to the third player and the points
at which its maxima are reached are summarized in Table IV below.%
\begin{equation*}
\begin{array}{c}
\text{Table IV: Payoffs to the third player and points of maxima.} \\ 
\begin{tabular}{cccccccc}
\hline
\multicolumn{2}{c}{} & \multicolumn{3}{c}{$l=8$} & \multicolumn{3}{c}{$l=5$}
\\ \hline
range of $q_{1}\&q_{2}$ & NE & $\text{payoff}$ & max & a$\text{t }\left(
q_{1},q_{2}\right) $ & payoff & max & a$\text{t }\left( q_{1},q_{2}\right) $
\\ \hline
$%
\begin{array}{c}
\frac{3q_{2}+6}{17}\leq q_{1} \\ 
\frac{3q_{1}+6}{17}\leq q_{2}%
\end{array}%
$ & AA & $4q_{1}+4q_{2}$ & $8$ & $\left( 1,1\right) $ & $q_{1}+q_{2}$ & $2$
& $\left( 1,1\right) $ \\ \hline
$%
\begin{array}{c}
q_{1}\leq \frac{3}{5} \\ 
q_{2}\leq \frac{3}{5}%
\end{array}%
$ & FF & $6q_{1}+6q_{2}$ & $\frac{108}{17}$ & $%
\begin{array}{c}
\left( \frac{3}{5},\frac{39}{85}\right)  \\ 
\left( \frac{39}{85},\frac{3}{5}\right) 
\end{array}%
$ & $3q_{1}+3q_{2}$ & $\frac{54}{17}$ & $%
\begin{array}{c}
\left( \frac{3}{5},\frac{39}{85}\right)  \\ 
\left( \frac{39}{85},\frac{3}{5}\right) 
\end{array}%
$ \\ \hline
$%
\begin{array}{c}
\frac{3}{5}\leq q_{1} \\ 
q_{2}\leq \frac{3q_{1}+6}{17}%
\end{array}%
$ & AF & $10q_{2}$ & $\frac{90}{17}$ & $\left( 1,\frac{9}{17}\right) $ & $%
7q_{2}-3q_{1}$ & $\frac{24}{17}$ & $\left( \frac{3}{5},\frac{39}{85}\right) $
\\ \hline
$%
\begin{array}{c}
q_{1}\leq \frac{3q_{2}+6}{17} \\ 
\frac{3}{5}\leq q_{2}%
\end{array}%
$ & FA & $10q_{1}$ & $\frac{90}{17}$ & $\left( \frac{9}{17},1\right) $ & $%
7q_{1}-3q_{2}$ & $\frac{24}{17}$ & $\left( \frac{39}{85},\frac{3}{5}\right) $
\\ \hline
\end{tabular}%
\end{array}%
\end{equation*}%
We can draw plots (see figure 2 and 3) that represent the payoff obtained by
the third player when he chooses different levels of attack. Comparing the
values at the points where the maxima are reached, the optimal choice of
attack levels is $q_{1}=1$ and $q_{2}=1$ if $l=8.$ Such a strategy of the
third player forces the first two players to form (or maintain) an alliance.
If $l=5,$ we see that the third player obtains maximum payoff at the points $%
\left( \frac{3}{5},\frac{39}{85}\right) $ or $\left( \frac{39}{85},\frac{3}{5%
}\right) $ but if the levels of attack are set to be exactly at these values
the payoff bi-matrix~(\ref{3Players_new_payoffs}) for the two firms is
non-generic and there exist an infinite number of Nash Equilibria for this
game. It includes the \textquotedblleft Ally\textquotedblright
-\textquotedblleft Ally\textquotedblright\ Nash Equilibrium which, if it is
played by the first two players, gives them the highest payoffs. The third
player obtains the payoff of $\frac{54}{17}$ only if the first two players
choose to fight. Therefore to prevent the first and the second players from
swapping between Nash Equilibria the third player must not choose $\left( 
\frac{3}{5},\frac{39}{85}\right) $ or $\left( \frac{39}{85},\frac{3}{5}%
\right) $ but choose $\left( \frac{3}{5}-\delta _{1},\frac{39}{85}-\delta
_{2}\right) $ or $\left( \frac{39}{85}-\delta _{1},\frac{3}{5}-\delta
_{2}\right) $ such that these points are still in (\textquotedblleft
Fight\textquotedblright -\textquotedblleft Fight\textquotedblright ) Nash
Equilibrium region and $\delta _{1}$ and $\delta _{2}$ are very small. We
can see that for this example the optimal choice of parameters $q_{1}$ and $%
q_{2}$ does not exist. The best strategy for the third player would be to
break the alliance by using appropriate levels of attack. In this case a
moderate, asymmetric attack will achieve higher payoff than a strong
(symmetric or asymmetric) attack.%
\begin{equation*}
\begin{array}{ccc}
\begin{tabular}{c}
{\unitlength=0.5400pt \begin{picture}(265.00,290.00)(0.00,0.00)
\put(140.00,260.00){\framebox(110.00,20.00){$AA$, $M$, $FF$}}
\put(140.00,260.00){\vector(0,-1){130.00}}
\put(200.00,80.00){\makebox(0.00,0.00){$AF$}}
\put(70.00,80.00){\makebox(0.00,0.00){$FF$}}
\put(190.00,190.00){\makebox(0.00,0.00){$AA$}}
\put(70.00,190.00){\makebox(0.00,0.00){$FA$}}
\put(260.00,30.00){\makebox(0.00,0.00){$q_1$}}
\put(30.00,260.00){\makebox(0.00,0.00){$q_2$}}
\put(230.00,10.00){\makebox(0.00,0.00){$1.0$}}
\put(10.00,230.00){\makebox(0.00,0.00){$1.0$}}
\put(10.00,190.00){\makebox(0.00,0.00){$0.8$}}
\put(10.00,150.00){\makebox(0.00,0.00){$0.6$}}
\put(10.00,110.00){\makebox(0.00,0.00){$0.4$}}
\put(10.00,70.00){\makebox(0.00,0.00){$0.2$}}
\put(190.00,10.00){\makebox(0.00,0.00){$0.8$}}
\put(150.00,10.00){\makebox(0.00,0.00){$0.6$}}
\put(110.00,10.00){\makebox(0.00,0.00){$0.4$}}
\put(70.00,10.00){\makebox(0.00,0.00){$0.2$}}
\put(110.00,110.00){\line(1,6){20.00}}
\put(110.00,110.00){\line(4,1){120.00}}
\put(230.00,230.00){\line(0,-1){200.00}}
\put(30.00,230.00){\line(1,0){200.00}}
\put(150.00,150.00){\line(0,-1){120.00}}
\put(30.00,150.00){\line(1,0){120.00}} \put(30.00,230.00){\line(-1,0){5.00}}
\put(30.00,190.00){\line(-1,0){5.00}} \put(30.00,150.00){\line(-1,0){5.00}}
\put(30.00,110.00){\line(-1,0){5.00}} \put(30.00,70.00){\line(-1,0){5.00}}
\put(230.00,30.00){\line(0,-1){5.00}} \put(190.00,30.00){\line(0,-1){5.00}}
\put(150.00,30.00){\line(0,-1){5.00}} \put(110.00,30.00){\line(0,-1){5.00}}
\put(70.00,30.00){\line(0,-1){5.00}} \put(30.00,30.00){\vector(1,0){220.00}}
\put(30.00,30.00){\vector(0,1){220.00}} \end{picture}} \\ 
$%
\begin{array}{c}
\text{Figure 1:} \\ 
\text{NE regions for player 1 and 2.}%
\end{array}%
$%
\end{tabular}
& 
\begin{tabular}{c}
{\unitlength=0.5400pt \begin{picture}(265.00,285.00)(0.00,0.00)
\put(140.00,260.00){\framebox(110.00,20.00){$4q_1+4q_2$}}
\put(140.00,260.00){\vector(0,-1){130.00}}
\put(200.00,80.00){\makebox(0.00,0.00){$10q_2$}}
\put(80.00,80.00){\makebox(0.00,0.00){$6q_1+6q_2$}}
\put(185.00,190.00){\makebox(0.00,0.00){$4q_1+4q_2$}}
\put(70.00,190.00){\makebox(0.00,0.00){$10q_1$}}
\put(260.00,30.00){\makebox(0.00,0.00){$q_2$}}
\put(30.00,260.00){\makebox(0.00,0.00){$q_2$}}
\put(230.00,10.00){\makebox(0.00,0.00){$1.0$}}
\put(10.00,230.00){\makebox(0.00,0.00){$1.0$}}
\put(10.00,190.00){\makebox(0.00,0.00){$0.8$}}
\put(10.00,150.00){\makebox(0.00,0.00){$0.6$}}
\put(10.00,110.00){\makebox(0.00,0.00){$0.4$}}
\put(10.00,70.00){\makebox(0.00,0.00){$0.2$}}
\put(190.00,10.00){\makebox(0.00,0.00){$0.8$}}
\put(150.00,10.00){\makebox(0.00,0.00){$0.6$}}
\put(110.00,10.00){\makebox(0.00,0.00){$0.4$}}
\put(70.00,10.00){\makebox(0.00,0.00){$0.2$}}
\put(110.00,110.00){\line(1,6){20.00}}
\put(110.00,110.00){\line(4,1){120.00}}
\put(230.00,230.00){\line(0,-1){200.00}}
\put(30.00,230.00){\line(1,0){200.00}}
\put(150.00,150.00){\line(0,-1){120.00}}
\put(30.00,150.00){\line(1,0){120.00}} \put(30.00,230.00){\line(-1,0){5.00}}
\put(30.00,190.00){\line(-1,0){5.00}} \put(30.00,150.00){\line(-1,0){5.00}}
\put(30.00,110.00){\line(-1,0){5.00}} \put(30.00,70.00){\line(-1,0){5.00}}
\put(230.00,30.00){\line(0,-1){5.00}} \put(190.00,30.00){\line(0,-1){5.00}}
\put(150.00,30.00){\line(0,-1){5.00}} \put(110.00,30.00){\line(0,-1){5.00}}
\put(70.00,30.00){\line(0,-1){5.00}} \put(30.00,30.00){\vector(1,0){220.00}}
\put(30.00,30.00){\vector(0,1){220.00}} \end{picture}} \\ 
$%
\begin{array}{c}
\text{Figure 2:} \\ 
\text{3 player's payoffs }l\text{=8.}%
\end{array}%
$%
\end{tabular}
& 
\begin{tabular}{c}
{\unitlength=0.5400pt \begin{picture}(265.00,285.00)(0.00,0.00)
\put(140.00,260.00){\framebox(110.00,20.00){$q_1+q_2$}}
\put(140.00,260.00){\vector(0,-1){130.00}}
\put(190.00,80.00){\makebox(0.00,0.00){$7q_2-3q_1$}}
\put(80.00,80.00){\makebox(0.00,0.00){$3q_1+3q_2$}}
\put(185.00,190.00){\makebox(0.00,0.00){$q_1+q_2$}}
\put(80.00,190.00){\makebox(0.00,0.00){$7q_1-3q_2$}}
\put(260.00,30.00){\makebox(0.00,0.00){$q_1$}}
\put(30.00,260.00){\makebox(0.00,0.00){$q_2$}}
\put(230.00,10.00){\makebox(0.00,0.00){$1.0$}}
\put(10.00,230.00){\makebox(0.00,0.00){$1.0$}}
\put(10.00,190.00){\makebox(0.00,0.00){$0.8$}}
\put(10.00,150.00){\makebox(0.00,0.00){$0.6$}}
\put(10.00,110.00){\makebox(0.00,0.00){$0.4$}}
\put(10.00,70.00){\makebox(0.00,0.00){$0.2$}}
\put(190.00,10.00){\makebox(0.00,0.00){$0.8$}}
\put(150.00,10.00){\makebox(0.00,0.00){$0.6$}}
\put(110.00,10.00){\makebox(0.00,0.00){$0.4$}}
\put(70.00,10.00){\makebox(0.00,0.00){$0.2$}}
\put(110.00,110.00){\line(1,6){20.00}}
\put(110.00,110.00){\line(4,1){120.00}}
\put(230.00,230.00){\line(0,-1){200.00}}
\put(30.00,230.00){\line(1,0){200.00}}
\put(150.00,150.00){\line(0,-1){120.00}}
\put(30.00,150.00){\line(1,0){120.00}} \put(30.00,230.00){\line(-1,0){5.00}}
\put(30.00,190.00){\line(-1,0){5.00}} \put(30.00,150.00){\line(-1,0){5.00}}
\put(30.00,110.00){\line(-1,0){5.00}} \put(30.00,70.00){\line(-1,0){5.00}}
\put(230.00,30.00){\line(0,-1){5.00}} \put(190.00,30.00){\line(0,-1){5.00}}
\put(150.00,30.00){\line(0,-1){5.00}} \put(110.00,30.00){\line(0,-1){5.00}}
\put(70.00,30.00){\line(0,-1){5.00}} \put(30.00,30.00){\vector(1,0){220.00}}
\put(30.00,30.00){\vector(0,1){220.00}} \end{picture}} \\ 
$%
\begin{array}{c}
\text{Figure 3:} \\ 
\text{3 player's payoffs }l\text{=5.}%
\end{array}%
$%
\end{tabular}%
\end{array}%
\end{equation*}

It is possible to show that (for $v=5,$ $x=2,$ $y=2,$ $f=2,$ $c=3,$ $%
\varepsilon =1,$ $\epsilon =1$) if the difference between value of the
market $l$ and the price of the attack $\delta $\ is

\begin{itemize}
\item less than $-1$ (note that in this case the actual cost of the attack
is $\sigma -\epsilon $ so the difference between value of the market $l$ and
the cost of the attack is less than $0)$ then it is optimal not to attack,
that is $q_{1}=0$ and $q_{2}=0;$

\item between $-1$ and $\frac{13}{4}$ then the optimal choice of parameters $%
q_{1}$ and $q_{2}$ does not exist and the third player should follow the
approach described above;

\item greater than $\frac{13}{4}$ then the optimal choice of attack levels
is $q_{1}=1$ and $q_{2}=1$; this choice of the third player forces the first
and the second players to form (or maintain) an alliance.
\end{itemize}

\noindent The proof can be obtained using item-by-item examination of all
possible cases. Similar conditions can be obtained for any generic set of
parameters.

\section{Summary and discussions.}

This simple example leads to a few important conclusions.

Firstly, as it can be seen from the results summarized in figure 1, the
successful formation of an alliance greatly depends on the position of the
players on the Nash Equilibrium diagram (figure 1). As in a single
interaction prisoners' dilemma in the absence of pressure cooperative
behaviour is not a Nash Equilibrium and alliances which are formed under
such conditions are not stable. However under the threat of an attack by
another player (or under unfavorable conditions) the players may be induced
to cooperate. It is clear that cooperative behaviour becomes a Nash
Equilibrium of the game~(\ref{3Players_new_payoffs}) because an attack on
the players changes their payoffs. But on the other hand, if, in modelling,
such changes are not taken into account we would obtain the result that
cooperation is not rational when, in fact, it is. In the considered model
strong pressure, where its level does not vary significantly between the
players will stabilize the possibility of forming an alliance. On the other
hand, an asymmetric, weak or moderate pressure does not provide an incentive
for cooperation.

In this instance the possibility for stabilizing an alliance arises by
applying a careful managing strategy. For example, assume the first player
position and suppose that $q_{1}=1/5$ (other values of $q_{1}$ can be
treated in a similar way).\ Assume also that the first player is interested
in forming an alliance. The idea is that the first player may be able to
change its own value of $q_{1}$ (for example by engaging in unilateral
commitments up to some extent)\ in such a way that the alliance becomes a
Nash Equilibrium. Three different approaches must be taken in such a
situation depending on the value on $q_{2}.$ If $0\leq q_{2}\leq 3/7$\ ($3/7$
appears as the point of intersection of $q_{1}=\left( 3q_{2}+6\right) /17$
and $q_{2}=\left( 3q_{1}+6\right) /17$) both players are in the
\textquotedblleft Fight\textquotedblright -\textquotedblleft
Fight\textquotedblright\ region and the pressure on the second player is not
high enough to promote a stable alliance. In this case the first player's
best strategy would be \textquotedblleft Fight\textquotedblright . The
unilateral commitments would only shift the players position into
\textquotedblleft Ally\textquotedblright -\textquotedblleft
Fight\textquotedblright\ region and the first company would suffer losses.
If $3/7<q_{2}\leq 3/5$ both players are still in the \textquotedblleft
Fight\textquotedblright -\textquotedblleft Fight\textquotedblright\ region
but the pressure on the second player is high enough to allow for the
possibility of an alliance. As can be seen from the figure~1, the first
company may change its own value of $q_{1}$ in such a way that
\textquotedblleft Ally\textquotedblright -\textquotedblleft
Ally\textquotedblright\ becomes a Nash Equilibrium. The unilateral
commitments in this case can lead to ensuring an alliance, though it should
be noted that the level of such commitments must be carefully assessed in
order to avoid shifting the players even further to the right on the Nash
Equilibrium diagram which would result in \textquotedblleft
Ally\textquotedblright -\textquotedblleft Fight\textquotedblright\
equilibrium and alliance failure. If $3/5\leq q_{2}$ the second player is
under high pressure and is likely to cooperate. The main threat to the
alliance is that the first player might prefer to use such an opportunity to
maximize its immediate payoff jeopardizing the future of the alliance. To
avoid such a contempt (for example the first company might be interested in
a strategic alliance which is expected to be highly beneficial in the
future) the first player might engage in unilateral commitments to stabilize
the alliance. As we can see engaging in unilateral commitments might be
highly beneficial for formation and maintaining of an alliance if used
appropriately. Note that it is assumed that the parameters of the model can
be estimated accurately enough. Although it is always possible to allow for
a margin of error, such an error must not be neglected when a strategy is
chosen as it might significantly change the strategy. The analysis of the
risks arising from this situation must be done, however we do not address
this issue in this paper.

The alliances must be frequently reassessed and managed$^{1}.$ In this
instance the results of the model can also be used for projections of
possible outcomes of an alliance. The estimates of future threats and
pressures must be made. If a company anticipates the position of an alliance
to leave the \textquotedblleft Ally\textquotedblright -\textquotedblleft
Ally\textquotedblright\ region it should be assessed if it is possible to
save the alliance by changing the values of $q_{1}$ or (and) $q_{2}.$ If it
is possible, careful management might save the alliance. If it is not the
case the company might save itself from heavy losses in the future by
withdrawing from an alliance prior to its failure.

Another result of the above analysis concerns the third player's strategy.
As we have seen from the examples, depending on the relationship between the
value of the market and the costs of the attack, the third player may either
wish to split the alliance or it may be better for the third player to
engage in the strongest symmetric attack. A strong attack, where the level
of attack does not vary significantly between the target firms, will
stabilize an alliance. An asymmetric, weak or moderate attack will result in
an alliance failure if the alliance is not managed well. It must be noticed
though that\ in this situation there exist mechanisms implying that good
management can always ensure that the alliance withstands an attack. The
strategy here is the same as when $3/5\leq q_{2}:$ the weakly attacked
player (let us say player one) must increase the value of $q_{1}$ to shift
the players position to the right on the Nash Equilibrium diagram. This
opportunity might be important then the third player represents hostile
conditions for one of the partners which are of a temporary nature. In this
situation the temporary sacrifice of another player might save the alliance
from failure bringing strategic benefits in the future.

As it can be seen from the example we have considered, a simple idea of
including environmental conditions into consideration provides several
important implications which would be overseen otherwise. In the considered
example, in addition to explaining the mechanisms which might be responsible
for formation of alliances, one of the important benefits is to allow us to
distinguish between two \textquotedblleft Fight\textquotedblright
-\textquotedblleft Fight\textquotedblright\ situations. In one unilateral
commitments are beneficial and in another would imply losses. These allow us
to use the proposed approach for determining successful strategies for
forming and maintaining alliances.

1. R. Gulati, T. Khanna and N. Nohria (1994) Unilateral commitments and the
importance of process in alliances, \textit{Sloan Management Rev.,} \textbf{%
35,} 61--69.

2. S.H. Park and G.R. Ungson (2001) Interfirm Rivalry and Managerial
Complexity: A Conceptual Framework of Alliance Failure, \textit{Organisation
Science,} \textbf{12} (1), 37-53.

3. E. Anderson (1990) Two Firms, One Frontier: On Assessing Joint Venture
Performance, \textit{Sloan Management Rev.,} \textbf{31} (2), 19-30.

4. J.M. Geringer and L. Hebert (1989) Control and Performance of
International Joint Ventures, \textit{J. of International Business Studies,} 
\textbf{20} (2), 235-254.

5. J.M. Geringer and L. Hebert (1989) Measuring Performance of International
Joint Ventures, \textit{J. of International Business Studies,} \textbf{22}
(2), 249-263.

6. K.R. Harrigan (1986) Managing for Joint Venture Success, Lexington, MA:
Lexington Books.

7. J.F. Hennart (1988) A Transition Cost Theory of Equity Joint Ventures, 
\textit{Strategic Management J.,} \textbf{9} (4), 361-374

8. A. Parkhe (1993) "Messy" Research, Methodological Predispositions, and
Theory Development in International Joint Ventures, \textit{Academy of
Management Rev.,} \textbf{18} (2), 227-268.

9. Y. Zhang and N. Rajagopalan (2002)\ Inter-partner Credible Threat in
International Joint Ventures: An Infinitely Repeated Prisoner's Dilemma
Model, \textit{J. of International Business Studies,} \textbf{33 }(3),
457-478.

10. J. Nash (1950) Equilibrium points in n-person games, \textit{Proc. Nat.
Acad. Sci.,} \textbf{36,} 48--49.

11. D. Fudenberg and J. Tirole (1991) Game theory, The MIT Press, Cambridge
Massachusetts USA.

12. J. von Neumann and O. Morgenstern (1944) Theory of games and economic
behaviour, Princeton Univ. Press, Princeton, NJ.

13. R. Axelrod and W. D. Hamilton (1981) The evolution of co-operation, 
\textit{Science,} \textbf{211,} 1390--1396.

14. M. Nowak and K. Sigmund (1993) A strategy of win-stay, lose-shift that
outperforms tit-for-tat in the prisoner's dilemma game, \textit{Nature,} 
\textbf{364,} 56--58.

15. J. Lorberaum (1994) No strategy is evolutionary stable, \textit{J. of
Theor. Biol.,} \textbf{168,} 117--130.

16. H. Gintis (2000) Game theory evolving: a problem-centered introduction
to modeling strategic interaction, Princeton Univ. Press, Princeton, New
Jersey.

17. R. Hoffmann (2001) The ecology of cooperation, \textit{Theory and
Decisions,} \textbf{50 }(2), 101-118.

18. D. P. Kraines and V. Y. Kraines (2000) Natural selection on memory-one
strategies for the iterated prisoner's dilemma, \textit{J. theor. Biol.,} 
\textbf{203,} 335-355.

19. J.Y. Wakano and N. Yamamura (2001) A simple learning strategy that
realizes robust cooperation better than Pavlov in iterated prisoners'
dilemma, \textit{J. of Ethology,} \textbf{19 }(1) 1-8.

20. M.A. Nowak and R.M. May (1992) Evolutionary games and spatial chaos, 
\textit{Nature,} \textbf{359} 826-829.

21. Y.F. Lim, K. Chen and C. Jayaprakash (2002) Scale invariant behaviour in
a spatial game of prisoners' dilemma, \textit{Physical Review E,} \textbf{65 
}(2), art. no. 026134,.

22. H.J. Jensen (1998) Self-organized criticality: emergent complex
behaviour in physical and biological systems, \textit{Cambridge lecture
notes in physics,} \textbf{10}, Cambridge Univ. Press.

23. M.N. Szilagyi and Z.C. Szilagyi (2002) Non-trivial solutions to the
N-person prisoners' dilemma, \textit{System Res\textit{earch }and
Behavioural Science,} \textbf{19 }(3), 281-290.

24. M. Kandori (2001) Introduction to repeated games with private
monitoring, \textit{J. of Economic Theory,} \textbf{102,} 1-15.

25. V. N. Kolokoltsov (1998) Nonexpansive maps and option pricing theory, 
\textit{Kybernetika,} \textbf{34} (6), 713--724.

\end{document}